\documentclass[aos,MSNbibl,seceqn,nameyear,dvips]{arximspdf}
\usepackage{graphicx}

\doi{10.1214/12-AOS1067} 
\volume{40}
\issue{6}
\pubyear{2012}
\firstpage{3137}
\lastpage{3160}

\makeatletter
\newcommand{\rrvert}{\vert}
\newcommand{\llvert}{\vert}

\newtheorem{theorem}{Theorem}
\newtheorem{cor}{Corollary}
\newtheorem{lem}{Lemma}
\newtheorem{prop}{Proposition}

\makeatother

\begin{document}
\begin{frontmatter}

\title{Improved multivariate normal mean estimation~with unknown
covariance when \lowercase{$p$} is greater than~\lowercase{$n$}}
\runtitle{Improved normal mean estimation for $p$ $>$ $n$}

\begin{aug}
\author[A]{\fnms{Didier} \snm{Ch\'etelat}\ead[label=e1]{dc623@cornell.edu}}
\and
\author[A]{\fnms{Martin T.} \snm{Wells}\corref{}\thanksref{t1}\ead[label=e2]{mtw1@cornell.edu}}
\runauthor{D. Ch\'etelat and M. T. Wells}
\affiliation{Cornell University}
\address[A]{Department of Statistical Science\\
Comstock Hall\\
Cornell University\\
Ithaca, New York 14853\\
USA\\
\printead{e1}\\
\hphantom{E-mail: }\printead*{e2}} 
\end{aug}

\thankstext{t1}{Supported in part by NSF Grant DMS-12-08488.}

\received{\smonth{1} \syear{2012}}
\revised{\smonth{11} \syear{2012}}

%
\begin{abstract}
We consider the problem of estimating the mean vector of a $p$-variate
normal $(\theta, \Sigma)$ distribution under invariant quadratic loss,
$(\delta-\theta)'\Sigma^{-1}(\delta-\theta)$, when the covariance is
unknown. We propose a new class of estimators that dominate the usual
estimator $\delta^0(X)=X$. The proposed estimators of $\theta$ depend
upon $X$ and an independent Wishart matrix $S$ with $n$ degrees of
freedom, however, $S$ is singular almost surely when $p > n$. The
proof of domination involves the development of some new unbiased
estimators of risk for the $p > n$ setting. We also find some
relationships between the amount of domination and the magnitudes of
$n$ and $p$.
\end{abstract}

%

\end{frontmatter}

\section{Introduction} \label{intro}

Suppose a $p$-dimensional random vector $X$ is observed which is
normally distributed, with mean vector $\theta$ and unknown positive
definite covariance matrix $\Sigma$, and we wish to estimate $\theta$
under the invariant quadratic loss
%
\begin{equation}
\label{eq2} L(\theta, \delta) = (\delta-\theta)'
\Sigma^{-1}(\delta-\theta).
\end{equation}
Since the covariance matrix $\Sigma$ is unknown, a random matrix $S$
is observed along with $X$, which is assumed to be independent of $X$,
and has a Wishart distribution with $n$ degrees of freedom, where $p >
n$. In high-dimensional estimation problems, where $p$, the number of
features, is nearly as
large as or larger than $n$, the number of observations, the ordinary
least squares estimator
does not typically provide a satisfactory estimate of $\theta$.

Modern data sets are increasingly becoming characterized by a number of
features that are much larger than the number of sample units
(large-$p$, small-$n$)
in contrast to classical data sets where the number of sample units is
often much larger than the number of random variables (small-$p$,
large-$n$). Modern applications in the $p > n$ setting include examples
from microarrays, association mapping, proteomics, radiology,
biomedical imaging, signal processing, climate modeling and finance.
For instance, in the case of microarray data,
the dimensionality is frequently in the thousands or beyond, while the
sample size is typically in the order of tens. The large-$p$, small-$n$
scenario poses challenges in most inferential settings. We are
considering a canonical setting. For the usual multivariate
location-scale estimation problem let $W = (W_1,\ldots, W_p)$ denote an
$N \times p$ matrix of data ($N$ is the number of observations and $p$
the number of features), where $W_i$ are taken from a $p$-dimensional
normal distribution with mean vector $\theta$ and covariance matrix
$\Xi$. In this article we let the $X$ and $S$ be the sample mean and
covariance of the features, respectively. In the context of this
notation, $\Sigma=N^{-1}\Xi$ and $n=N-1$.

The usual estimator under invariant quadratic loss is $\delta_0(X)=X$.
It is minimax and admissible when $p \leq2$ and $p \leq n$. However,
when $p\ge3$ and $p \leq n$, $\delta_0(X)$ remains minimax but is no
longer admissible. Explicit improvements are known in the multivariate
normal case [\citet{BerBock76Ann,Beretal77Ann,BerHaff83SD},
Gleser (\citeyear{Gleser79Ann,Gleser86Ann}), \citet{JamesStein61proBer}]
and in the case of elliptically
symmetric distribution [\citet{SriBil89JMA,fourdrinier03}].

In this article we primarily concentrate on the case $p > n $ and
construct a class of estimators, depending on the sufficient statistics
$(X,S)$, of the form
%
\begin{equation}\label{eq3}
\delta(X,S)=X+g(X,S),
\end{equation}
which dominate $\delta_0(X)$ under invariant quadratic loss. Note that,
although the loss in (\ref{eq2}) is invariant, the estimate in
(\ref{eq3}) may not be [except for $\delta_0(X)$]. This class
generalizes several estimators studied previously for the multivariate
normal distribution to the $p \leq n$ setting
[\citet{BerBock76Ann,Beretal77Ann,BerHaff83SD}, Gleser
(\citeyear{Gleser79Ann,Gleser86Ann}), \citet{JamesStein61proBer}].
Examples of estimators we study here in this setting extend the class
of so-called Baranchik estimators and includes a new high-dimensional
James--Stein estimator
\[
\delta^{\mathrm{JS}}_a(X,S)= \biggl(I-\frac{aSS^+}{X'S^+X} \biggr)X,
\]
where $0\leq a \leq\frac{2(n-2)}{p-n+3}$ and $S^+$ is the
Moore--Penrose inverse of $S$.

The estimation of the inverse covariance matrix, namely, the precision
matrix $\Sigma^{-1}$,
of a multivariate normal distribution has been an important problem in
practical situations
as well as from a theoretical perspective. But, when $p > n$,
the Wishart-distributed sample covariance matrix is singular; in this
case, one is tempted\vadjust{\goodbreak} to construct estimators using
the Moore--Penrose generalized inverse $S^+$.
Recently there has been an increased interest in the problem of
estimating the covariance matrix of large dimension given variables of
dimension larger
than the number of observations [\citet
{Bickel08,Asp08,konno09,Ledoit,Levina08,Rothman08}].\looseness=-1

Our method of proof relies on an unbiased estimator of risk difference, say,
$\rho(X,S)$. Specifically, we show that, for $g(X,S)$ of the form
$-\frac{r(X'S^+X)SS^+}{X'S^+X}X$, the estimator
$\delta(X,S)=X+g(X,S)$ dominates $X$ provided $\rho(X,S) \le0 $. In
the next section we present the main results and their proofs are given
in Section~\ref{sectech}. We need Stein's integration-by-parts
identity [\citet{Stein81Ann}] and the so-called Stein--Haff
identity for the singular Wishart distribution. The Stein--Haff identity
was derived by \citet{Haff79JMA} and \citet{Stein77} for
the full rank Wishart distribution. A~similar identity for the elliptically contoured model has been given by
\citet{fourdrinier03}. We make some concluding comments in
Section~\ref{secnumstud}.

For a matrix $M$, let $M'$ denote its transpose, $M^+$ its
Moore--Penrose pseudo-inverse and $\frac{\partial M}{\partial t}$ its
componentwise derivative matrix, that is, the matrix such that $
(\frac{\partial M}{\partial t} )_{ij}=\frac{\partial
M_{ij}}{\partial t}$. Moreover, let $\delta_{ij}$ denote the Kronecker
delta.\looseness=-1

\section{Main results}
Let $X$ be a random vector distributed as $N_p(\theta,\Sigma)$ with
unknown $\theta$ and $\Sigma$.
Suppose an estimator of $\Sigma$ is available, say, $S \sim\mathrm
{Wishart}_p(n,\Sigma)$, with $S$ independent of $ X$. By definition of
the Wishart distribution, we can write $S=Y'Y$ for some matrix normal
$Y \sim N_{n\times p}(0,I\otimes\Sigma)$. An elementary property of
this distribution is that $S$ is (almost surely) invertible if $p \leq
n$, and (almost surely) singular if $p > n$ [cf. \citet{srivastava79}].

An usual estimator of $\theta$ is $\delta^0(X,S)=X$; however, it
turns out that this estimator is inadmissible under quadratic loss. If
some estimator $S\sim\mathrm{Wishart}_p(n,\Sigma)$ is available,
with $n\geq p\geq3$, $\delta^0$ is dominated by the so-called
James--Stein estimator
\[
\delta^{\mathrm{JS}}(X,S)= \biggl(1-\frac{(p-2)/(n-p+3)}{X'S^{-1}X} \biggr)X.
\]
The main contribution of this article is to extend this type of result
to a more general class of estimators in the $p > n$ setting.

For some positive, bounded and differentiable function $r\dvtx  \mathbb{R}
\rightarrow\mathbb{R}$, define the Baranchik-type estimator
%
\begin{eqnarray}
\label{deltar} \delta_r(X,S)&=& \biggl(I-\frac{r(X'S^+X)SS^+}{X'S^+X} \biggr)X
\nonumber\\[-9pt]\\[-9pt]
&=&X+g(X,S),
\nonumber
\end{eqnarray}
where $I$ is the identity matrix and $S^+$ denotes the Moore--Penrose
inverse of $S$. This estimator generalizes the usual \citet
{Baranchik} estimator to the unknown covariance setting for
$p>n$.\vadjust{\goodbreak}

\begin{theorem}\label{domination} Let $\min(p, n)\geq3$. Suppose that:
\begin{longlist}[(iii)]
\item[(i)] \label{dom1} $r$ satisfies $0\leq r \leq\frac{2 (\min
(n, p)-2)}{n+p-2 \min(n, p)+3}$;
\item[(ii)] $r$ is nondecreasing; and
\item[(iii)] $r'$ is bounded.
\end{longlist}
Then under invariant quadratic loss, $\delta_r$ dominates $\delta^0$.
\end{theorem}
%
Throughout the article we will use the expression $\operatorname{tr}(SS^+)$, which of
course equals $\min(n, p)$. This notation allows us to simultaneously
handle both the $p>n$ and $n\geq p$ cases. The condition $\min
(p,n)\geq3$ merely guarantees that condition (i) of Theorem \ref
{domination} holds for some $r$ and is reminiscent of the dimension
cutoff in classical Stein estimation.
\begin{pf*}{Proof of Theorem~\ref{domination}}
The hypotheses of the theorem imply that $r$ is differentiable almost
everywhere. Under invariant quadratic loss, the difference in risk
between $\delta_r$ and $\delta^0$ is given by
%
\begin{eqnarray}\label{thmdel}
\Delta_\theta &=& E_\theta \bigl[\bigl(X+g(X,S)-\theta
\bigr)'\Sigma^{-1}\bigl(X+g(X,S)-\theta\bigr) \bigr]
\nonumber
\\
&&{}-E_\theta \bigl[(X-\theta)'\Sigma^{-1}(X-\theta)
\bigr]
\\
&=&2E_\theta \bigl[g(X,S)'\Sigma^{-1}(X-\theta)
\bigr]+E_\theta \bigl[g(X,S)'\Sigma^{-1}g(X,S)
\bigr].\nonumber
\end{eqnarray}

In order to show the domination result, we need to show that under the
sufficient conditions on $r$, (\ref{thmdel}) is nonpositive for all
$\theta$. First, for the leftmost term of (\ref{thmdel}) it can be
shown that
\[
2E_\theta \bigl[g(X,S)'\Sigma^{-1}(X-\theta)
\bigr] = 2E_\theta \bigl[ \mathrm{div}_X g(X,S) \bigr].
\]
\citet{fourdrinier03} give a more general form of this result in
their Lemma 1(i); it is essentially an extension of Stein's classical
integration by parts identity. By using Lemma~\ref{divx} in Section~\ref{sectech},
we have that
%
\begin{eqnarray}\label{thmmaina}\quad
2E_\theta \bigl[ \mathrm{div}_X g(X,S) \bigr] &=&
-2E_\theta \biggl[ \mathrm{div}_X \frac{r(X'S^+X)SS^+X}{X'S^+X} \biggr]
\nonumber\\[-8pt]\\[-8pt]
&=&-2E_\theta \biggl[ 2r'\bigl(X'S^+X\bigr) +
r\bigl(X'S^+X\bigr)\frac{\operatorname
{tr}(SS^+)-2}{X'S^+X} \biggr].\nonumber
\end{eqnarray}

For the right term of (\ref{thmdel}), we find, through Lemma \ref
{divkonno} in Section~\ref{sectech},
\begin{eqnarray*}
&&
E_\theta \bigl[g(X,S)'\Sigma^{-1}g(X,S) \bigr]
\\
&&\qquad= E_\theta \biggl[\operatorname{tr} \biggl(\Sigma^{-1}S
r^2\bigl(X'S^+X\bigr) \frac{S^+XX'S^+S}{(X'S^+X)^2} \biggr) \biggr]
\\
&&\qquad= E_\theta\biggl[ n \operatorname{tr} \biggl(r^2
\bigl(X'S^+X\bigr)\frac
{S^+XX'S^+S}{(X'S^+X)^2} \biggr)
\\
&&\qquad\quad{}+ \operatorname{tr} \biggl(Y'\nabla_Y \biggl\{
r^2\bigl(X'S^{+}X\bigr)\frac{SS^{+}XX'S^{+}}{(X'S^{+}X)^2}
\biggr\} \biggr) \biggr].
\end{eqnarray*}
The finiteness of the risk of $\delta_r$ is guaranteed to hold by
Theorem~\ref{momthm} in Section~\ref{sectech} for all $p$ and $n$.

Now applying Lemma~\ref{trdelta} in Section~\ref{sectech}, we find
%
\begin{eqnarray}\label{thmmainb}
&&E_\theta\biggl[ n \operatorname{tr} \biggl(r^2
\bigl(X'S^{+}X\bigr)\frac
{S^+XX'S^+S}{(X'S^+X)^2} \biggr)
\nonumber
\\
&&\hspace*{5.6pt}\quad{}+ \operatorname{tr} \biggl(Y'\nabla_Y \biggl\{
r^2\bigl(X'S^{+}X\bigr)\frac{SS^{+}XX'S^{+}}{(X'S^{+}X)^2}
\biggr\} \biggr)\biggr]
\nonumber
\\
&&\hspace*{5.6pt}\qquad= E_\theta\biggl[ n \frac{r^2(X'S^{+}X)}{X'S^+X} -4r\bigl(X'S^{+}X
\bigr)r'\bigl(X'S^{+}X\bigr)
\\
&&\hspace*{47pt}\hspace*{5.6pt}\qquad\quad{} +r^2\bigl(X'S^{+}X\bigr)\frac{p-2 \operatorname
{tr}(SS^{+})+3}{X'S^{+}X}
\biggr]
\nonumber
\\
&&\hspace*{5.6pt}\qquad= E_\theta \biggl[r^2\bigl(X'S^{+}X
\bigr)\frac{n+p-2 \operatorname
{tr}(SS^{+})+3}{X'S^{+}X}-4r\bigl(X'S^{+}X
\bigr)r'\bigl(X'S^{+}X\bigr) \biggr].\nonumber
\end{eqnarray}
Replacing (\ref{thmmaina}) and (\ref{thmmainb}) back into (\ref
{thmdel}), we obtain
\begin{eqnarray*}
\Delta_\theta &=& E_\theta\biggl[ r^2
\bigl(X'S^{+}X\bigr)\frac{n+p-2
\operatorname{tr}(SS^{+})+3}{X'S^{+}X}
\\
&&\hspace*{17pt}{}-2r\bigl(X'S^+X\bigr)\frac{\operatorname{tr}(SS^+)-2}{X'S^+X}
\\
&&\hspace*{30pt}{} -4 r'\bigl(X'S^+X\bigr)
\bigl\{ 1+r\bigl(X'S^{+}X\bigr) \bigr\}\biggr].
\end{eqnarray*}
Since $r$ is nonnegative and nondecreasing, it follows that $-4
r'(X'S^+X) \{1+r(X'S^{+}X)\}\leq0$.
Finally, for the $X$ and $S$ such that $r(X'S^+X) \neq0$,
\begin{eqnarray*}
&& r^2\bigl(X'S^{+}X\bigr)\frac{n+p-2 \operatorname{tr}(SS^{+})+3}{X'S^{+}X}
-2r\bigl(X'S^+X\bigr)\frac{\operatorname{tr}(SS^+)-2}{X'S^+X} \leq0
\\
&&\quad\Leftrightarrow\quad r\bigl(X'S^{+}X\bigr) \leq
\frac{2 (\operatorname
{tr}(SS^+)-2)}{n+p-2 \operatorname{tr}(SS^+)+3}= \frac{2 (\min
(n,p)-2)}{n+p-2 \min(n,p)+3}.
\end{eqnarray*}
Therefore, under the three sufficient conditions on $r$, it follows
that $\Delta_\theta\leq0$ for any~$\theta$, that is,
the domination result holds.
\end{pf*}

In the $p>n$ setting, we obtain the following two corollaries.
%
\begin{cor}\label{coro1} For $p>n\geq3$, $\delta_r$ dominates $\delta^0$ under
invariant quadratic loss for all $r$ nondecreasing, differentiable and
satisfying
%
\begin{equation}\label{cor1}
0\leq r \leq\frac{2(n-2)}{p-n+3}.
\end{equation}
\end{cor}
%
\begin{cor}[(James--Stein estimator with large $p$ and small $n$)] For
$p>n\geq3$ and $a \in\mathbb{R}$, the James--Stein-like estimator
%
\begin{equation}\label{cor2}
\delta^{\mathrm{JS}}_a(X,S)= \biggl(I-\frac{aSS^+}{X'S^+X} \biggr)X
\end{equation}
dominates $\delta^0$ under invariant quadratic loss for all
\[
0\leq a \leq\frac{2(n-2)}{p-n+3}.
\]
\end{cor}

Note that if $p$ is only moderately larger than $n$, Corollary~\ref{coro1}
implies that one can construct an estimator with substantial
improvement over $\delta^0$. However, in the ultra-high-dimensional
setting the denominator in (\ref{cor1}) could be quite large and,
consequently, the amount of improvement over $\delta^0$ could be quite
small. The estimator in (\ref{cor2}) generalizes the classical
James--Stein with unknown covariance matrix,
\[
\delta^{\mathrm{JS}}_a(X,S)= \biggl(1-\frac{a}{X'S^{-1}X} \biggr)X,
\]
which is, of course, restricted to the case $p \leq n$, for $a\in
\mathbb{R}_+$. In this setting, this result is consistent with
previous bounds in \citet{fourdrinier03} (where $n-1$ is used
instead of our $n$).

\section{Technical results and proofs}\label{sectech}
It remains to clarify several of the somewhat technical computations
used in the proof of Theorem~\ref{domination}. We provide them in this section; these
computations are likely to be of independent interest and showcase
several technical maneuvers that the reader could find useful in
dealing with singular Wishart matrices.
%
\begin{prop}\label{compprop} Let $Y$ be an $n\times p$ matrix,
$S=Y'Y$, $X$ a $p$ vector and $F=X'S^+X$. It then follows that
\begin{eqnarray*}
&&\mbox{\textup{\hphantom{ii}(i)}}\quad 
\biggl\{\frac{\partial S}{\partial Y_{\alpha\beta
}} \biggr\}_{kl}=\delta_{\beta k} Y_{\alpha l} + \delta_{\beta
l} Y_{\alpha k};
\\
&&\mbox{\textup{\hphantom{i}(ii)}}\quad 
\frac{\partial F}{\partial Y_{\alpha\beta
}}=-2\bigl(X'S^{+}Y'\bigr)_\alpha\bigl(S^{+}X\bigr)_\beta+ 2\bigl(X'S^{+}S^{+}Y'\bigr)_\alpha
\bigl(\bigl(I-SS^{+}\bigr)X\bigr)_\beta;
\\
&&\mbox{\textup{(iii)}}\quad 
\frac{\partial \{ S^{+}XX'SS^{+} \}_{kl}}{\partial Y_{\alpha\beta}}
\\
&&\hphantom{\textup{(iii)}}\qquad
=\bigl(S^{+}S^{+}Y'\bigr)_{k \alpha} \bigl(
\bigl(I-SS^{+}\bigr)XX'SS^{+}
\bigr)_{\beta l}
\\
&&\hphantom{\textup{(iii)}}\qquad\quad
{}-S^{+}_{k\beta}\bigl(YS^{+}XX'SS^{+}
\bigr)_{\alpha l}-\bigl(S^{+}Y'\bigr)_{k\alpha
}
\bigl(S^{+}XX'SS^{+}\bigr)_{\beta l}
\\
&&\hphantom{\textup{(iii)}}\qquad\quad
{}+\bigl(I-SS^{+}\bigr)_{k \beta} \bigl(YS^{+}S^{+}XX'SS^{+}
\bigr)_{\alpha l}
\\
&&\hphantom{\textup{(iii)}}\qquad\quad
{}+\bigl(S^{+}XX'\bigr)_{k \beta}
\bigl(YS^{+}\bigr)_{\alpha l} + \bigl(S^{+}XX'Y'
\bigr)_{k \alpha
}\bigl(S^{+}\bigr)_{\beta l}
\\
&&\hphantom{\textup{(iii)}}\qquad\quad
{}+\bigl(S^{+}XX'S^{+}Y'
\bigr)_{k\alpha}\bigl(I-SS^{+}\bigr)_{\beta l}
\\
&&\hphantom{\textup{(iii)}}\qquad\quad
{}-\bigl(S^{+}XX'SS^{+}\bigr)_{k\beta}
\bigl(YS^{+}\bigr)_{\alpha l} -\bigl(S^{+}XX'SS^{+}Y'
\bigr)_{k\alpha}\bigl(S^{+}\bigr)_{\beta l}.
\end{eqnarray*}
\end{prop}
\begin{pf}
First, notice that from the usual chain-rule that
\[
\biggl\{\frac{\partial S}{\partial Y_{\alpha\beta}} \biggr\}_{kl} =
\frac{\partial}{\partial Y_{\alpha\beta}}S_{kl} = \frac{\partial
}{\partial Y_{\alpha\beta}} \sum
_q Y_{qk}Y_{ql} = \delta_{\beta
k}Y_{\alpha l}
+ \delta_{\beta l}Y_{\alpha k}.
\]
This shows (i).

Let $A$ be a symmetric matrix and $t \in\mathbb{R}$, then
\begin{eqnarray*}
\frac{\partial A^{+}}{\partial t}&=& -A^{+}\,\frac{\partial A}{\partial t}A^{+} +
\bigl(I-AA^{+}\bigr)\,\frac{\partial A}{\partial t}A^{+}A^{+}\\
&&{}+A^{+}A^{+}\,\frac{\partial A}{\partial t}\bigl(I-AA^{+}\bigr).
\end{eqnarray*}
This result was, it seems, first proved in \citet{golub73}, as
their Theorem 4.3, but can be found in standard textbooks on elementary
linear algebra. Also, again for $A$ symmetric, we have $AA^{+}=A^{+}A$
and $ A(I-AA^{+})=(I-AA^{+})A=A^{+}(I-AA^{+})=(I-AA^{+})A^{+}=0$. This
easily follows from elementary properties of the Moore--Penrose pseudoinverse.

Since $S=Y'Y$, notice through a singular value decomposition argument
that $SS^+Y'=Y'$ and, thus, $(I-SS^+)Y'=0$. Using (i), we find that
\begin{eqnarray*}
\frac{\partial F}{\partial Y_{\alpha\beta}} &=& X' \,\frac{\partial
S^{+}}{\partial Y_{\alpha\beta}} X
\\
&=& -\sum_{k,l} \bigl(X'S^{+}
\bigr)_k \{\delta_{\beta k} Y_{\alpha l} +
\delta_{\beta l}Y_{\alpha k} \}\bigl(S^{+}X
\bigr)_l
\\
&&{}+ \sum_{k,l} \bigl(X'S^{+}S^{+}
\bigr)_k \{\delta_{\beta k} Y_{\alpha l} +
\delta_{\beta l} Y_{\alpha k} \} \bigl(\bigl(I-SS^{+}
\bigr)X\bigr)_l
\\
&&{}+ \sum_{k,l} \bigl(X'
\bigl(I-SS^{+}\bigr)\bigr)_k \{\delta_{\beta k}
Y_{\alpha l} + \delta_{\beta l} Y_{\alpha k} \}
\bigl(S^{+}S^{+}X\bigr)_l
\\
&=&-\sum_{l} \bigl(X'S^{+}
\bigr)_\beta Y_{\alpha l} \bigl(S^{+}X
\bigr)_l - \sum_{k}
\bigl(X'S^{+}\bigr)_k Y_{\alpha k}
\bigl(S^{+}X\bigr)_\beta
\\
&&{}+\sum_{l} \bigl(X'S^{+}S^{+}
\bigr)_\beta Y_{\alpha l} \bigl(\bigl(I-SS^{+}\bigr)X
\bigr)_l
\\
&&{}+ \sum_{k} \bigl(X'S^{+}S^{+}
\bigr)_k Y_{\alpha k} \bigl(\bigl(I-SS^{+}\bigr)X
\bigr)_\beta
\\
&&{}+\sum_{l} \bigl(X'
\bigl(I-SS^{+}\bigr)\bigr)_\beta Y_{\alpha l}
\bigl(S^{+}S^{+}X\bigr)_l
\\[-1pt]
&&{}+ \sum_{k} \bigl(X'
\bigl(I-SS^{+}\bigr)\bigr)_k Y_{\alpha k}
\bigl(S^{+}S^{+}X\bigr)_\beta
\\[-1pt]
&=&-2\bigl(X'S^{+}Y'\bigr)_\alpha
\bigl(S^{+}X\bigr)_\beta+2\bigl(X'S^{+}S^{+}Y'
\bigr)_\alpha \bigl(\bigl(I-SS^{+}\bigr)X\bigr)_\beta,
\end{eqnarray*}
which gives (ii).

Using (i), we have that for any conformable matrices $A$ and $B$
\begin{eqnarray*}
\biggl( A \,\frac{\partial S}{\partial Y_{\alpha\beta}} B \biggr)_{kl} &=&\sum
_{i,j} A_{ki} \biggl\{\frac{\partial S}{\partial
Y_{\alpha\beta}} \biggr
\}_{ij} B_{jl}
\\[-1pt]
&=&\sum_{i,j} A_{ki} \{
\delta_{\beta i} Y_{\alpha j} + \delta_{\beta j}Y_{\alpha i}
\} B_{jl}
\\[-1pt]
&=&\sum_{j} A_{k\beta} Y_{\alpha j}
B_{jl} + \sum_{i} A_{ki}
Y_{\alpha i} B_{\beta l}
\\[-1pt]
&=&A_{k \beta} (YB)_{\alpha l} + \bigl(AY'
\bigr)_{k \alpha} B_{\beta l}.
\end{eqnarray*}
Therefore, using again $(I-SS^+)Y'=0$,
\begin{eqnarray*}
&&
\frac{\partial\{ S^{+}XX'SS^{+}\}_{kl}}{\partial
Y_{\alpha\beta}} \\[-1pt]
&&\qquad=\biggl\{ S^{+}S^{+}\,
\frac{\partial S}{\partial Y_{\alpha\beta
}}\bigl(I-SS^{+}\bigr)XX'SS^{+}
\\[-1pt]
&&\hspace*{3.6pt}\qquad\quad{}-S^{+}\,\frac{\partial S}{\partial Y_{\alpha\beta}}S^{+}XX'SS^{+}
+\bigl(I-SS^{+}\bigr)\,\frac{\partial S}{\partial Y_{\alpha\beta
}}S^{+}S^{+}XX'SS^{+}
\\[-1pt]
&&\hspace*{3.6pt}\qquad\quad{}+ S^{+}XX'\,\frac{\partial S}{\partial Y_{\alpha\beta}}S^{+} +
S^{+}XX'SS^{+}S^{+}\,\frac{\partial S}{\partial Y_{\alpha\beta
}}
\bigl(I-SS^{+}\bigr)
\\[-1pt]
&&\hspace*{3.6pt}\qquad\quad{}-S^{+}XX'SS^{+}\,\frac{\partial S}{\partial Y_{\alpha\beta}}S^{+}
+S^{+}XX'S\bigl(I-SS^{+}\bigr)
\,\frac{\partial S}{\partial Y_{\alpha\beta
}}S^{+}S^{+} \biggr\}_{kl}
\\[-1pt]
&&\qquad= \bigl(S^{+}S^{+}Y'
\bigr)_{k \alpha} \bigl(\bigl(I-SS^{+}\bigr)XX'SS^{+}
\bigr)_{\beta l}
\\[-1pt]
&&\qquad\quad{}-S^{+}_{k\beta}\bigl(YS^{+}XX'SS^{+}
\bigr)_{\alpha l}-\bigl(S^{+}Y'\bigr)_{k\alpha
}
\bigl(S^{+}XX'SS^{+}\bigr)_{\beta l}
\\[-1pt]
&&\qquad\quad{}+\bigl(I-SS^{+}\bigr)_{k \beta} \bigl(YS^{+}S^{+}XX'SS^{+}
\bigr)_{\alpha l}
\\[-1pt]
&&\qquad\quad{}+\bigl(S^{+}XX'\bigr)_{k \beta}
\bigl(YS^{+}\bigr)_{\alpha l} + \bigl(S^{+}XX'Y'
\bigr)_{k \alpha
}\bigl(S^{+}\bigr)_{\beta l}
\\[-1pt]
&&\qquad\quad{}+\bigl(S^{+}XX'S^{+}Y'
\bigr)_{k\alpha}\bigl(I-SS^{+}\bigr)_{\beta l}
\\[-1pt]
&&\qquad\quad{}-\bigl(S^{+}XX'SS^{+}\bigr)_{k\beta}
\bigl(YS^{+}\bigr)_{\alpha l} -\bigl(S^{+}XX'SS^{+}Y'
\bigr)_{k\alpha}\bigl(S^{+}\bigr)_{\beta l},
\end{eqnarray*}
which gives (iii).\vadjust{\goodbreak}
\end{pf}
%
\begin{lem}\label{trdelta}
Under the hypotheses of Theorem~\ref{domination}, we have
\begin{eqnarray*}
&&\operatorname{tr} \biggl(Y'\nabla_Y \biggl\{
r^2\bigl(X'S^{+}X\bigr)\frac
{SS^{+}XX'S^{+}}{(X'S^{+}X)^2}
\biggr\} \biggr)
\\
&&\qquad=-4r\bigl(X'S^{+}X\bigr)r'
\bigl(X'S^{+}X\bigr)+r^2\bigl(X'S^{+}X
\bigr)\frac{p-2\operatorname{tr}(SS^{+})+3}{X'S^{+}X},
\end{eqnarray*}
where $\nabla_Y$ is interpreted as the matrix with components $(\nabla_Y)_{ij}=\frac{\partial}{\partial Y_{ij}}$.
\end{lem}
\begin{pf}
To simplify computations, in what will follows, we let $F \equiv
X'S^{+}X$. We then have
%
\begin{eqnarray}
&& \biggl[Y'\nabla_Y \biggl\{ r^2(F)
\frac
{SS^{+}XX'S^{+}}{F^2} \biggr\} \biggr]_{ij}
\nonumber
\\
&&\qquad=\sum_{\alpha,\beta}\bigl(Y'
\bigr)_{i\alpha}\,\frac{\partial}{\partial
Y_{\alpha\beta}} \biggl\{ r^2(F)
\frac{(SS^{+}XX'S^{+})_{\beta
j}}{F^2} \biggr\}
\nonumber
\\
\label{lemmaina}
&&\qquad=2 \sum_{\alpha,\beta}\bigl(Y'
\bigr)_{i\alpha} r(F)r'(F) \,\frac{\partial
F}{\partial Y_{\alpha\beta}} \cdot
\frac{(SS^{+}XX'S^{+})_{\beta j}}{F^2}
\\
\label{lemmainb}
&&\qquad\quad{}+\sum_{\alpha,\beta}\bigl(Y'
\bigr)_{i\alpha} r^2(F)\frac{({\partial
}/{\partial Y_{\alpha\beta}})\{(SS^{+}XX'S^{+})_{\beta
j}\}}{F^2}
\\
\label{lemmainc}
&&\qquad\quad{}+\sum_{\alpha,\beta}\bigl(Y'
\bigr)_{i\alpha} r^2(F)\frac{- 2\,({\partial
F}/{\partial Y_{\alpha\beta}})
(SS^{+}XX'S^{+})_{\beta j}}{F^3}.
\end{eqnarray}
To simplify (\ref{lemmaina}) and (\ref{lemmainc}), we apply
Proposition~\ref{compprop}(ii) to get
\begin{eqnarray*}
&&\sum_{\alpha,\beta}\bigl(Y'
\bigr)_{i\alpha} \biggl\{\frac{\partial
F}{\partial Y_{\alpha\beta}} \biggr\}
\bigl(SS^{+}XX'S^{+}\bigr)_{\beta
j}
\\
&&\qquad=-2\sum_{\alpha,\beta}\bigl(Y'
\bigr)_{i\alpha}\bigl(X'S^{+}Y'
\bigr)_\alpha \bigl(S^{+}X\bigr)_\beta
\bigl(SS^{+}XX'S^{+}\bigr)_{\beta j}
\\
&&\qquad\quad{}+2\sum_{\alpha,\beta} \bigl(X'S^{+}S^{+}Y'
\bigr)_\alpha(Y)_{\alpha i} \bigl(S^{+}XX'SS^{+}
\bigr)_{j\beta} \bigl(\bigl(I-SS^{+}\bigr)X\bigr)_\beta
\\
&&\qquad=-2X'S^{+}X \bigl(SS^{+}XX'S^{+}
\bigr)_{ij}.
\end{eqnarray*}
Using this, we get for (\ref{lemmaina})
%
\begin{eqnarray}\label{lemsuba}
&&
2 \sum_{\alpha,\beta}\bigl(Y'
\bigr)_{i\alpha} r(F)r'(F) \,\frac{\partial
F}{\partial Y_{\alpha\beta}} \cdot
\frac{(SS^{+}XX'S^{+})_{\beta
j}}{F^2}
\nonumber\\[-8pt]\\[-8pt]
&&\qquad=-4r(F)r'(F) \frac{(SS^{+}XX'S^{+})_{ij}}{F}\nonumber
\end{eqnarray}
and (\ref{lemmainc}) becomes
%
\begin{eqnarray}\label{lemsubc}
&&\sum_{\alpha,\beta}\bigl(Y'\bigr)_{i\alpha}
r^2(F) \frac{- 2\,({\partial
F}/{\partial Y_{\alpha\beta}}) \cdot(SS^{+}XX'S^{+})_{\beta j}}{F^3}
\nonumber\\[-8pt]\\[-8pt]
&&\qquad= 4r^2(F) \frac{(SS^{+}XX'S^{+})_{ij}}{F^2}.\nonumber
\end{eqnarray}

This leaves the term (\ref{lemmainb}) to analyze. Using Proposition
\ref{compprop}(iii),
\begin{eqnarray*}
&&
\sum_{\alpha,\beta}\bigl(Y'
\bigr)_{i\alpha}\,\frac{\partial}{\partial
Y_{\alpha\beta}} \bigl\{\bigl(SS^{+}XX'S^{+}
\bigr)_{\beta j} \bigr\}
\\
&&\qquad=\sum_{\alpha,\beta}\bigl(Y'
\bigr)_{i\alpha}\,\frac{\partial\{
S^{+}XX'SS^{+}\}_{j\beta}}{\partial Y_{\alpha\beta}}
\\
&&\qquad=\sum_{\alpha,\beta} \bigl\{
\bigl(S^{+}S^{+}Y'\bigr)_{j \alpha
}Y_{\alpha i}
\bigl(\bigl(I-SS^{+}\bigr)XX'SS^{+}
\bigr)_{\beta\beta}
\\
&&\hspace*{17pt}\qquad\quad{}-S^{+}_{j \beta}\bigl(Y'\bigr)_{i\alpha}
\bigl(YS^{+}XX'SS^{+}\bigr)_{\alpha\beta}
\\
&&\hspace*{17pt}\qquad\quad{}-\bigl(S^{+}Y'\bigr)_{j \alpha}Y_{\alpha i}
\bigl(S^{+}XX'SS^{+}\bigr)_{\beta\beta}
\\
&&\hspace*{17pt}\qquad\quad{}+\bigl(I-SS^{+}\bigr)_{j \beta} \bigl(Y'
\bigr)_{i\alpha}\bigl(YS^{+}S^{+}XX'SS^{+}
\bigr)_{\alpha
\beta}
\\
&&\hspace*{17pt}\qquad\quad{}+\bigl(S^{+}XX'\bigr)_{j \beta}
\bigl(Y'\bigr)_{i\alpha}\bigl(YS^{+}
\bigr)_{\alpha\beta}
\\
&&\hspace*{17pt}\qquad\quad{}+ \bigl(S^{+}XX'Y'\bigr)_{j \alpha}Y_{\alpha i}
\bigl(S^{+}\bigr)_{\beta\beta}
\\
&&\hspace*{17pt}\qquad\quad{}+\bigl(S^{+}XX'S^{+}Y'
\bigr)_{j \alpha}Y_{\alpha i} \bigl(I-SS^{+}
\bigr)_{\beta\beta}
\\
&&\hspace*{17pt}\qquad\quad{}-\bigl(S^{+}XX'SS^{+}\bigr)_{j \beta}
\bigl(Y'\bigr)_{i\alpha}\bigl(YS^{+}
\bigr)_{\alpha\beta}
\\
&&\hspace*{65pt}\qquad\quad{}-\bigl(S^{+}XX'SS^{+}Y'
\bigr)_{j \alpha}Y_{\alpha i} \bigl(S^{+}\bigr)_{\beta\beta}
\bigr\}
\\
&&\qquad= \bigl(S^{+}XX'SS^{+}
\bigl(I-SS^{+}\bigr)\bigr)_{ij}
\\
&&\qquad\quad{}-\bigl(SS^{+}XX'S^{+}\bigr)_{ij}
-\operatorname{tr}\bigl(S^{+}XX'SS^{+}\bigr)
\bigl(SS^+\bigr)_{ij}
\\
&&\qquad\quad{}+\operatorname{tr}\bigl(\bigl(I-SS^{+}\bigr)XX'SS^{+}
\bigr) \bigl(S^{+}\bigr)_{ij}
\\
&&\qquad\quad{}+\bigl(SS^{+}XX'S^{+}\bigr)_{ij}
+\operatorname{tr}\bigl(S^{+}\bigr) \bigl(SXX'S^{+}
\bigr)_{ij}
\\
&&\qquad\quad{}+\operatorname{tr}\bigl(I-SS^{+}\bigr) \bigl(SS^{+}XX'S^{+}
\bigr)_{ij}
\\
&&\qquad\quad{}-\bigl(SS^{+}XX'S^{+}\bigr)_{ij }
-\operatorname{tr}\bigl(S^{+}\bigr) \bigl(SXX'S^{+}
\bigr)_{ij}
\\
&&\qquad= \bigl(p-\operatorname{tr}\bigl(SS^{+}\bigr)-1\bigr)
\bigl\{ SS^{+}XX'S^{+} \bigr
\}_{ij}-\bigl(X'S^{+}X\bigr) \bigl\{
SS^{+} \bigr\}_{ij}.
\end{eqnarray*}

Next, applying this computation in (\ref{lemmainb}), we obtain
%
\begin{eqnarray}\label{lemsubb}\qquad
&&
\sum_{\alpha,\beta}\bigl(Y'\bigr)_{i\alpha}
r^2(F) \frac{({\partial
}/{\partial Y_{\alpha\beta}})\{(SS^{+}XX'S^{+})_{\beta j}\}
}{F^2}
\nonumber
\\[-8pt]\\[-8pt]
&&\qquad= \bigl(p-\operatorname{tr}\bigl(SS^{+}\bigr)-1\bigr)
r^2(F)\frac{(
SS^{+}XX'S^{+})_{ij}}{F^2}
 -r^2(F) \frac{(SS^+)_{ij}}{F}.\nonumber
\end{eqnarray}

Now we can combine (\ref{lemsuba}), (\ref{lemsubb}) and (\ref
{lemsubc}) together to complete the proof. That is, we have
\begin{eqnarray*}
&&\operatorname{tr} \biggl(Y'\nabla_Y \biggl\{
r^2(F)\frac
{SS^{+}XX'S^{+}}{F^2} \biggr\} \biggr)
\\[-2pt]
&&\qquad=\sum_i \biggl\{-4r(F)r'(F)
\frac{(SS^{+}XX'S^{+})_{ii}}{F}
\\[-2pt]
&&\hspace*{18.1pt}\qquad\quad{}+ 4r^2(F) \frac{(SS^{+}XX'S^{+})_{ii}}{F^2}
\\[-2pt]
&&\hspace*{18.1pt}\qquad\quad{}+ \bigl(p-\operatorname{tr}\bigl(SS^{+}\bigr)-1\bigr)
r^2(F)\frac{(
SS^{+}XX'S^{+})_{ii}}{F^2}
-r^2(F) \frac{(SS^+)_{ii}}{F} \biggr\}
\\[-2pt]
&&\qquad=-4r(F)r'(F)+r^2(F)\frac{p-2 \operatorname{tr}(SS^{+})+3}{F}
\end{eqnarray*}
as desired.
\end{pf}
%
\begin{lem} \label{divx} Under the hypotheses of Theorem \ref
{domination} we have
\[
\mathrm{div}_X \frac{r(X'S^+X)SS^+X}{X'S^+X}=2r'
\bigl(X'S^+X\bigr) + r\bigl(X'S^+X\bigr)
\frac{\operatorname{tr}(SS^+)-2}{X'S^+X}.
\]
\end{lem}
\begin{pf}
Again, to simplify computations, let us denote $X'S^+X$ by $F$. We find
\begin{eqnarray*}
&&
\mathrm{div}_X  \biggl\{ r(F)\frac{SS^+X}{F} \biggr\}\\[-2pt]
&&\qquad=\sum_i \frac{\partial}{\partial X_i} \biggl\{ r(F)
\frac
{(SS^+X)_i}{F} \biggr\}
\\[-2pt]
&&\qquad=\sum_i r'(F) \,\frac{\partial F}{\partial X_i}
\frac{(SS^+X)_i}{F}
+ r(F)\frac{({\partial}/{\partial X_i}) \{
(SS^+X)_i \}}{F} \\[-2pt]
&&\qquad\quad{}- r(F)\frac{({\partial F}/{\partial X_i}) (SS^+X)_i}{F^2}
\\[-2pt]
&&\qquad=\sum_i r'(F) \biggl\{
\frac{\partial}{\partial X_i} \sum_{k,l} X_k
X_l S^+_{kl} \biggr\}\frac{(SS^+X)_i}{F}
\\[-2pt]
&&\qquad\quad{}+ r(F)\frac{({\partial}/{\partial X_i}) \sum_k
(SS^+)_{ik}X_k}{F}
\\[-2pt]
&&\qquad\quad{}- r(F)\frac{ \{({\partial}/{\partial X_i} )\sum_{k,l}
X_k X_l S^+_{kl}  \}(SS^+X)_i}{F^2}
\\[-2pt]
&&\qquad=\sum_i r'(F) \bigl\{
\bigl(X'S^+\bigr)_i+\bigl(X'S^+
\bigr)_i \bigr\}\frac
{(SS^+X)_i}{F}
\\[-2pt]
&&\qquad\quad{}+ r(F)\frac{(SS^+)_{ii}}{F} - r(F)\frac{ \{(X'S^+)_i+(X'S^+)_i  \}\cdot
(SS^+X)_i}{F^2}
\\[-2pt]
&&\qquad= 2r'(F) + r(F)\frac{\operatorname{tr}(SS^+)-2}{F}
\end{eqnarray*}
as desired.
\end{pf}

The following result is an extension of a result in \citet
{konno09}. This type of result was first obtained by \citet
{KubSriv08} and then was extended by \citet{konno09}. In our
generalization we make use of a divergence version of Stein's lemma
that comes with somewhat weaker moment conditions, rather than the
element-by-element assumptions in \citet{konno09}. These weaker
moment conditions allow us to cover the $p$ equals $n$ and $n+1$ cases.
%
\begin{lem}\label{divkonno} Let $Y\sim N_{n\times p}(0,I_n\otimes
\Sigma)$, let $S=Y'Y$ which has, by definition,
a~$\mathrm{Wishart}_p(n,\Sigma)$ distribution, and let $G(S)$ be a
$p\times p$ random matrix that
depends on $S$. Let $\nabla_Y$ be interpreted as the matrix with components
$(\nabla_Y)_{ij}=\frac{\partial}{\partial Y_{ij}}$, and for $A$ the
symmetric positive definite square root of $\Sigma$,
define $\tilde{Y}=YA^{-1}$ and $H=AGA^{-1}$. Then
\[
E \bigl[\operatorname{tr} \bigl(\Sigma^{-1}SG \bigr) \bigr]=E \bigl[n
\operatorname{tr} (G )+\operatorname{tr} \bigl(Y'
\nabla_Y G' \bigr) \bigr]
\]
under the conditions
%
\begin{equation}\label{divkonnocond}
E \bigl[\bigl\llvert \mathrm{div}_{\operatorname{vec}(\tilde Y)} \cdot \operatorname{vec} (
\tilde Y H ) \bigr\rrvert \bigr]<\infty,
\end{equation}
where $\operatorname{vec}(M)$ denotes the vectorization of a matrix
$M$.
\end{lem}
\begin{pf}
Define $\tilde S=\tilde Y'\tilde Y=A^{-1}SA^{-1}$. Notice that, by
construction, $\tilde Y \sim N_{n\times p}(0,I_n\otimes I_p)$---this
means, by definition of the matrix normal distribution, that
$\operatorname{vec}(\tilde Y) \sim N_{np}(0,I_{np})$. We can write
\begin{eqnarray*}
E \bigl[\operatorname{tr} (\tilde SH ) \bigr]&=&E \biggl[\sum
_{\alpha,i,j}\tilde Y_{\alpha i} \tilde Y_{\alpha j}
H_{ji} \biggr]\\
&=&E \bigl[\operatorname{vec}(\tilde Y) \cdot
\operatorname{vec} (\tilde Y H ) \bigr].
\end{eqnarray*}
Using the divergence form of Stein's lemma, which can be found in Lemma~A.1 in\vadjust{\goodbreak} \citet{fourdrinier03b}, we obtain, under the moment
conditions outlined in (\ref{divkonnocond}),
\begin{eqnarray*}
E \bigl[\operatorname{vec}(\tilde Y) \cdot\operatorname{vec} (\tilde Y H )
\bigr] &=& E \bigl[\mathrm{div}_{\operatorname{vec}(\tilde Y)} \operatorname {vec} (\tilde Y H )
\bigr]
\\
&=& E \biggl[\sum_{\alpha,i,j}\frac{\partial}{\partial\tilde
Y_{\alpha i}} \tilde
Y_{\alpha j} H_{ji} \biggr]
\\
&=& E \biggl[\sum_{\alpha,i,j} \delta_{ij}
H_{ji} + \tilde Y_{\alpha j} \,\frac{\partial H_{ji}}{\partial\tilde Y_{\alpha
i}} \biggr]
\\
&=& E \biggl[n \sum_{i} H_{ii} + \sum
_{\alpha,i,j} \tilde Y_{\alpha j} \,\frac{\partial}{\partial\tilde Y_{\alpha i}}
H_{ji} \biggr].
\end{eqnarray*}
This last expression can be expressed in a compact matrix form as
\[
E \bigl[\operatorname{tr} (\tilde SH ) \bigr]=E \bigl[n \operatorname{tr}(H)+
\operatorname{tr} \bigl(\bigl(\tilde Y'\nabla_{\tilde Y}
\bigr)'H \bigr) \bigr].
\]
Finally, we notice
\begin{eqnarray*}
E \bigl[\operatorname{tr} (H ) \bigr]&=&E \bigl[\operatorname{tr}
\bigl(AGA^{-1} \bigr) \bigr],
\\
E \bigl[\operatorname{tr} (\tilde SH ) \bigr]&=&E \bigl[\operatorname{tr}
\bigl(A^{-1}SGA^{-1} \bigr) \bigr],
\\
E \bigl[\operatorname{tr} \bigl(\bigl(\tilde Y'
\nabla_{\tilde Y}\bigr)'H \bigr) \bigr]&=&E \bigl[\operatorname{tr}
\bigl(A\bigl(Y'\nabla_Y\bigr)'GA^{-1}
\bigr) \bigr],
\end{eqnarray*}
which concludes the proof.
\end{pf}
%
\begin{theorem}\label{momthm} Let $Y\sim N_{n\times p}(0,I_n\otimes
\Sigma
)$ and
for $A$ the symmetric positive definite square root of $\Sigma$, let
$\tilde{Y}=YA^{-1}$. Let $r$ be any bounded differentiable
nonnegative function $r\dvtx \mathbb{R} \rightarrow[0,C_1]$ with bounded
derivative $|r'|\leq C_2$. Define
\[
G=r^2\bigl(X'S^{+}X\bigr)\frac{S^{+}XX'S^{+}S}{(X'S^{+}X)^2}
\]
and $H=AGA^{-1}$. Then for all $p$ and $n$
%
\begin{equation}\label{momthmcond}
E \bigl[\bigl\llvert \mathrm{div}_{\operatorname{vec}(\tilde Y)} \operatorname{vec} (\tilde Y H
) \bigr\rrvert \bigr]<\infty.
\end{equation}
\end{theorem}
\begin{pf}
We first compute $\mathrm{div}_{\operatorname{vec}(\tilde Y)}
\operatorname{vec} (\tilde Y H )$. As always, to ease
notation, we shall write $F=X'S^{+}X$. We have
%
\begin{eqnarray}
&&
\mathrm{div}_{\operatorname{vec}(\tilde Y)} \operatorname {vec} (\tilde Y H ) \nonumber\\[-2pt]
&&\qquad=\sum
_{\alpha,i,j}\,\frac{\partial}{\partial\tilde Y_{\alpha i}} \{\tilde Y_{\alpha j}
H_{ji} \}
\nonumber
\\[-2pt]
&&\qquad=n\sum_{i}H_{ii}+\sum
_{\alpha,j}\tilde Y_{\alpha j} \,\frac{\partial
H_{ji}}{\partial\tilde Y_{\alpha i}}
\nonumber
\\[-2pt]
&&\qquad=n\sum_{i}H_{ii} +\sum
_{\alpha,\beta,i,j}\tilde Y_{\alpha j} A_{\beta i}
\,\frac{\partial}{\partial Y_{\alpha\beta}} \biggl\{ r^2(F)\frac{ \{ AS^{+}XX'SS^{+}A^{-1}
\}_{ji}}{F^2} \biggr
\}
\nonumber
\\[-2pt]
&&\qquad=n\sum_{i}H_{ii}+\sum
_{\alpha,\beta,i,j}\tilde Y_{\alpha j} A_{\beta i}
\nonumber
\\[-2pt]
\label{mom1a}
&&\hspace*{79pt}\qquad\quad{}\times\biggl\{2r(F)r'(F) \,\frac{\partial F}{\partial Y_{\alpha\beta
}} \frac{ \{ AS^{+}XX'SS^{+}A^{-1} \}_{ji}}{F^2}
\\[-2pt]
\label{mom1b}
&&\hspace*{96.3pt}\qquad\quad{}+\frac{r^2(F)}{F^2} \sum_{k,l} A_{jk}
\,\frac{\partial\{
S^{+}XX'SS^{+}\}_{kl}}{\partial Y_{\alpha\beta}} A^{-1}_{li}
\\[-2pt]
\label{mom1c}
&&\hspace*{96.3pt}\qquad\quad{} -r^2(F) \bigl\{ AS^{+}XX'SS^{+}A^{-1}
\bigr\}_{ji} \frac{2\,{\partial F}/{\partial Y_{\alpha\beta}}}{F^3} \biggr\}.
\end{eqnarray}
We simplify each part of the expression. For (\ref{mom1a}), using
Proposition~\ref{compprop}(ii), we find
%
\begin{eqnarray}\label{mom2a}
&&2 \sum_{\alpha,\beta,i,j}\tilde Y_{\alpha j}
A_{\beta i} r(F)r'(F) \,\frac{\partial F}{\partial Y_{\alpha\beta}} \frac{ \{
AS^{+}XX'SS^{+}A^{-1} \}_{ji}}{F^2}
\nonumber\hspace*{-15pt}
\\[-2pt]
&&\quad= 4\frac{r(F)r'(F)}{F^2}
\nonumber\hspace*{-15pt}
\\[-2pt]
&&\qquad{}\times\sum_{\alpha,\beta,i,j} \bigl\{-\bigl(X'S^{+}Y'\bigr)_\alpha
\tilde Y_{\alpha j} \bigl\{ AS^{+}XX'SS^{+}A^{-1}
\bigr\}_{ji} A_{i \beta}\bigl(S^{+}X
\bigr)_\beta
\nonumber\hspace*{-15pt}\\[-2pt]
&&\hspace*{41.1pt}\qquad{}+ \bigl(X'S^{+}S^{+}Y'
\bigr)_\alpha\tilde Y_{\alpha j} \bigl\{ AS^{+}XX'SS^{+}A^{-1}
\bigr\}_{ji} A_{i \beta
}\bigl(\bigl(I-SS^{+}\bigr)X
\bigr)_\beta \bigr\}\hspace*{-15pt}
\\[-2pt]
&&\quad=-4\frac{r(F)r'(F)}{F^2} \bigl(X'S^{+}Y'YA^{-1}AS^{+}XX'SS^{+}A^{-1}AS^{+}X
\bigr)
\nonumber\hspace*{-15pt}
\\[-2pt]
&&\qquad{}+4\frac{r(F)r'(F)}{F^2} \bigl(X'S^{+}S^{+}Y'YA^{-1}AS^{+}XX'SS^{+}A^{-1}A
\bigl(I-SS^{+}\bigr)X \bigr)
\nonumber\hspace*{-15pt}
\\[-2pt]
&&\quad=-4r(F)r'(F).\nonumber\hspace*{-15pt}
\end{eqnarray}
Similarly, for (\ref{mom1c})
%
\begin{eqnarray}\label{mom2c}
&&
\sum_{\alpha,\beta,i,j}\tilde Y_{\alpha j} A_{\beta i}
r^2(F) \bigl\{ AS^{+}XX'SS^{+}A^{-1}
\bigr\}_{ji} \frac{2\,
{\partial F}/{\partial Y_{\alpha\beta}}}{F^3}
\nonumber
\\[-2pt]
&&\qquad= 4\frac{r^2(F)}{F^3} \sum_{\alpha,\beta,i,j}
\bigl(X'S^{+}Y'\bigr)_\alpha \tilde
Y_{\alpha j} \bigl\{ AS^{+}XX'SS^{+}A^{-1}
\bigr\}_{ji}A_{i \beta}\bigl(S^{+}X
\bigr)_\beta
\nonumber\\[-8pt]\\[-8pt]
&&\qquad= 4\frac{r^2(F)}{F^3} \bigl(X'S^{+}Y'YA^{-1}AS^{+}XX'SS^{+}A^{-1}AS^{+}X
\bigr)
\nonumber
\\[-2pt]
&&\qquad= 4\frac{r^2(F)}{F}.\nonumber
\end{eqnarray}
This leaves us with (\ref{mom1b}). Using Proposition
\ref{compprop}(iii), we obtain
%
\begin{eqnarray}\label{mom2b}
&&
\sum_{\alpha,\beta,i,j}\tilde Y_{\alpha j} A_{\beta i}
\frac
{r^2(F)}{F^2} \sum_{k,l} A_{jk}\,
\frac{\partial \{
S^{+}XX'SS^{+} \}_{kl}}{\partial Y_{\alpha\beta}} A^{-1}_{li}
\nonumber
\\
&&\qquad=\frac{r^2(F)}{F^2} \sum_{\alpha,\beta,i,j,k,l} \tilde
Y_{\alpha
j} A_{\beta i} A_{jk} A^{-1}_{li}
\nonumber
\\
&&\hspace*{66.6pt}\qquad\quad{}\times\bigl\{\bigl(S^{+}S^{+}Y\bigr)_{k \alpha}
\bigl(\bigl(I-SS^{+}\bigr)XX'SS^{+}
\bigr)_{\beta l}
\nonumber
\\
&&\hspace*{83.3pt}\qquad\quad{}-S^{+}_{k\beta}\bigl(YS^{+}XX'SS^{+}
\bigr)_{\alpha l}
\nonumber
\\
&&\hspace*{83.3pt}\qquad\quad{}-\bigl(S^{+}Y\bigr)_{k\alpha}\bigl(S^{+}XX'SS^{+}
\bigr)_{\beta l}
\nonumber
\\
&&\hspace*{83.3pt}\qquad\quad{}+\bigl(I-SS^{+}\bigr)_{k \beta} \bigl(YS^{+}S^{+}XX'SS^{+}
\bigr)_{\alpha l}
\nonumber
\\
&&\hspace*{83.3pt}\qquad\quad{}+\bigl(S^{+}XX'\bigr)_{k \beta}
\bigl(YS^{+}\bigr)_{\alpha l}
\nonumber
\\
&&\hspace*{83.3pt}\qquad\quad{}+ \bigl(S^{+}XX'Y'\bigr)_{k \alpha}
\bigl(S^{+}\bigr)_{\beta l}
\nonumber
\\
&&\hspace*{83.3pt}\qquad\quad{}+\bigl(S^{+}XX'S^{+}Y'
\bigr)_{k\alpha}\bigl(I-SS^{+}\bigr)_{\beta l}
\nonumber
\\
&&\hspace*{83.3pt}\qquad\quad{}-\bigl(S^{+}XX'SS^{+}\bigr)_{k\beta}
\bigl(YS^{+}\bigr)_{\alpha l}
\nonumber
\\
&&\hspace*{120pt}\qquad\quad{}-\bigl(S^{+}XX'SS^{+}Y'
\bigr)_{k\alpha}\bigl(S^{+}\bigr)_{\beta l} \bigr\}
\nonumber
\\
&&\qquad=\frac{r^2(F)}{F^2} \sum_{\alpha,\beta,i,j,k,l} \bigl
\{
A_{jk}\bigl(S^{+}S^{+}Y\bigr)_{k \alpha}
\tilde Y_{\alpha j} A_{i\beta
}\bigl(\bigl(I-SS^{+}
\bigr)XX'SS^{+}\bigr)_{\beta l}A^{-1}_{li}
\nonumber
\\
&&\hspace*{82pt}\qquad{}-\tilde Y'_{j \alpha}\bigl(YS^{+}XX'SS^{+}
\bigr)_{\alpha l}A^{-1}_{li} A_{i
\beta}S^{+}_{\beta k}A_{kj}
\nonumber
\\
&&\hspace*{82pt}\qquad{}-A_{jk}\bigl(S^{+}Y\bigr)_{k\alpha}\tilde
Y_{\alpha j} A_{i\beta
}\bigl(S^{+}XX'SS^{+}
\bigr)_{\beta l}A^{-1}_{li}
\nonumber
\\
&&\hspace*{82pt}\qquad{}+ \tilde Y'_{j \alpha}\bigl(YS^{+}S^{+}XX'SS^{+}
\bigr)_{\alpha l}A^{-1}_{li} A_{i\beta}
\bigl(I-SS^{+}\bigr)_{\beta k}A_{kj}
\\
&&\hspace*{82pt}\qquad{}+ \tilde Y'_{j \alpha}\bigl(YS^{+}
\bigr)_{\alpha l}A^{-1}_{li} A_{i\beta
}
\bigl(XX'S^{+}\bigr)_{\beta k}A_{kj}
\nonumber
\\
&&\hspace*{82pt}\qquad{}+ A_{jk}\bigl(S^{+}XX'Y'
\bigr)_{k \alpha}\tilde Y_{\alpha j} A_{i\beta
}\bigl(S^{+}
\bigr)_{\beta l}A^{-1}_{li}
\nonumber
\\
&&\hspace*{82pt}\qquad{}+A_{jk}\bigl(S^{+}XX'S^{+}Y'
\bigr)_{k\alpha}\tilde Y_{\alpha j} A_{i\beta
}\bigl(I-SS^{+}
\bigr)_{\beta l}A^{-1}_{li}
\nonumber
\\
&&\hspace*{82pt}\qquad{}- \tilde Y'_{j \alpha}\bigl(YS^{+}
\bigr)_{\alpha l}A^{-1}_{li}A_{i\beta
}
\bigl(SS^{+}XX'S^{+}\bigr)_{\beta k}A_{kj}
\nonumber
\\
&&\hspace*{122pt}\qquad{}-A_{jk}\bigl(S^{+}XX'SS^{+}Y'
\bigr)_{k\alpha}\tilde Y_{\alpha j} A_{i\beta
}\bigl(S^{+}
\bigr)_{\beta l}A^{-1}_{li} \bigr\}
\nonumber
\\
&&\qquad=\frac{r^2(F)}{F^2} \bigl\{\operatorname {tr}
\bigl(AS^{+}S^{+}Y'YA^{-1}\bigr)\cdot
\operatorname {tr}\bigl(A\bigl(I-SS^{+}\bigr)XX'SS^{+}A^{-1}
\bigr)
\nonumber
\\
&&\hspace*{43.5pt}\qquad{}-\operatorname{tr}\bigl(A^{-1}Y'YS^{+}XX'SS^{+}A^{-1}AS^{+}A
\bigr)
\nonumber
\\
&&\hspace*{43.5pt}\qquad{}-\operatorname{tr}\bigl(AS^{+}Y'YA^{-1}\bigr)
\operatorname {tr}\bigl(AS^{+}XX'SS^{+}A^{-1}
\bigr)
\nonumber
\\
&&\hspace*{43.5pt}\qquad{}+\operatorname{tr}\bigl(A^{-1}Y'YS^{+}S^{+}XX'SS^{+}A^{-1}A
\bigl(I-SS^{+}\bigr)A\bigr)
\nonumber
\\
&&\hspace*{43.5pt}\qquad{}+\operatorname{tr}\bigl(A^{-1}Y'YS^{+}A^{-1}AXX'S^{+}A
\bigr)
\nonumber
\\
&&\hspace*{43.5pt}\qquad{}+\operatorname{tr}\bigl(AS^{+}XX'Y'YA^{-1}
\bigr)\cdot\operatorname {tr}\bigl(AS^{+}A^{-1}\bigr)
\nonumber
\\
&&\hspace*{43.5pt}\qquad{}+\operatorname{tr}\bigl(AS^{+}XX'S^{+}Y'YA^{-1}
\bigr) \operatorname {tr}\bigl(A\bigl(I-SS^{+}\bigr)A^{-1}
\bigr)
\nonumber
\\
&&\hspace*{43.5pt}\qquad{}-\operatorname{tr}\bigl(A^{-1}Y'YS^{+}A^{-1}ASS^{+}XX'S^{+}A
\bigr)
\nonumber
\\
&&\hspace*{85.3pt}\qquad{}-\operatorname{tr}\bigl(AS^{+}XX'SS^{+}Y'YA^{-1}
\bigr) \operatorname {tr}\bigl(AS^{+}A^{-1}\bigr) \bigr\}
\nonumber
\\
&&\qquad=\frac{r^2(F)}{F^2} \cdot\bigl\{-X'S^{+}X-
\operatorname {tr}\bigl(SS^{+}\bigr)\cdot X'S^{+}X
\nonumber
\\
&&\hspace*{50.3pt}\qquad{}+X'S^{+}X+X'SS^+X\cdot\operatorname{tr}
\bigl(S^+\bigr)+X'S^{+}X\cdot \bigl(p-\operatorname{tr}
\bigl(SS^+\bigr)\bigr)
\nonumber
\\
&&\hspace*{167.3pt}\qquad{}-X'S^+X-X'SS^+X\cdot\operatorname{tr}\bigl(S^+\bigr)
\bigr\}
\nonumber
\\
&&\qquad=\frac{r^2(F)}{F} \bigl(p- \operatorname{tr}\bigl(SS^{+}
\bigr)-1 \bigr).\nonumber
\end{eqnarray}
Having re-expressed $\mathrm{div}_{\operatorname{vec}(\tilde Y)}
\operatorname{vec} (\tilde Y H )$, we now need to bound it
above. By virtue of (\ref{mom2a}), (\ref{mom2c}) and (\ref{mom2b}),
we have
%
\begin{eqnarray}\label{mom3}
&&E\bigl[\bigl\llvert \mathrm{div}_{\operatorname{vec}(\tilde Y)} \operatorname{vec} (\tilde Y H
) \bigr\rrvert \bigr]
\nonumber
\\
&&\qquad= E\biggl[\biggl| n \operatorname{tr}(H)+4\frac{r^2(F)}{F}
\nonumber\\[-8pt]\\[-8pt]
&&\hspace*{15.1pt}\qquad\quad{} + \bigl(p-\operatorname{tr}\bigl(SS^{+}\bigr)-1 \bigr)
\frac
{r^2(F)}{F} -4r(F)r'(F) \biggr|\biggr]
\nonumber
\\
&&\qquad\leq C^2_1\bigl\llvert 3+p-\operatorname{tr}
\bigl(SS^{+}\bigr)+n\bigr\rrvert E \biggl[ \frac1{F} \biggr] +
4C_1C_2.\nonumber
\end{eqnarray}
It only remains to show that $E [ \frac1{F}  ]$ is finite.
By definition of the Wishart matrix distribution,
we can define a $T\sim\mathrm{Wishart}_p(n,I_n)$ such that $S=ATA$.
Let $T=H'DH$ be the spectral decomposition of $T$, with
$D=\operatorname{diag}(\lambda_i)$. Write the eigenvalues of $T^+$ as
$\lambda^+_i$, so that $D^{-1}=\operatorname{diag}(\lambda^+_i)$,
and let $\lambda^+_{\min}$ be the smallest nonzero eigenvalue of
$T^+$. The following two identities follow from \citet{tian04}
[Theorem 1.1, equations (1.2) and (1.4)] and symmetry of $T$:
\begin{eqnarray*}
(ATA)^+&=&\bigl(T^+TA\bigr)^+T^+\bigl(AT^+T\bigr)^+,
\\
\bigl(T^+TA\bigr)^+\bigl(T^+T\bigr)&=&\bigl(T^+TA\bigr)^+.
\end{eqnarray*}
Using these identities, we have
\begin{eqnarray*}
X'S^+X &=& X'(ATA)^+X=X'\bigl(T^+TA
\bigr)^+T^+\bigl(AT^+T\bigr)^+X
\\
&=& \sum_{k} \bigl\{ X'
\bigl(T^+TA\bigr)^+H' \bigr\}_k^2
\lambda^+_k
\\
&\geq& \lambda^+_{\min} \cdot X'\bigl(T^+TA
\bigr)^+H'H\bigl(AT^+T\bigr)^+X
\\
&=& \lambda^+_{\min} \cdot X'\bigl(T^+TA\bigr)^+\bigl(T^+T
\bigr) \bigl(AT^+T\bigr)^+X
\\
&=& \lambda^+_{\min} \cdot X'\bigl(T^+TA\bigr)^+\bigl(AT^+T
\bigr)^+X.
\end{eqnarray*}
Applying Cauchy--Schwarz provides us with the bound
\[
X'\bigl(T^+TA\bigr)^+\bigl(T^+TA\bigr)X \leq X'
\bigl(T^+TA\bigr)^+\bigl(AT^+T\bigr)^+ X X'\bigl(AT^+T\bigr)
\bigl(T^+TA\bigr)X
\]
so that we then have
\begin{eqnarray*}
\frac1F=\frac1{X'S^+X} &\leq& \frac1{\lambda^+_{\min}} \frac
1{X'\bigl(T^+TA\bigr)^+\bigl(AT^+T\bigr)^+X}
\\
&\leq& \frac1{\lambda^+_{\min}} \frac{X'AT^+TAX}{X'(T^+TA)^+(T^+TA)X}.
\end{eqnarray*}
To ease notation, let us write $Q=AT^+TA $ and $R=(T^+TA)^+(T^+TA)$.
Collecting the results together, we bound (\ref{mom3}) by
%
\begin{equation}\label{mom4}
\leq C^2_1\bigl\llvert 3+p-2 \operatorname{tr}
\bigl(SS^{+}\bigr)+n\bigr\rrvert E \biggl[ \frac1{\lambda^+_{\min}}
\frac{X'QX}{X'RX} \biggr] + 4C_1C_2.
\end{equation}
We now use some independence results. We can write the singular value
decomposition of $T$ as $T=H'DH$, but we can also write it as $T=H'_1D_1H_1$,
where $H_1$ is semi-orthogonal ($H_1H'_1=I)$ and $D_1$ is the matrix of
the positive eigenvalues of~$T$. If $T$ has full rank
(i.e., $n\geq p$), then this coincides with the singular value
decomposition of $T$.
In the full rank case, \citet{srivastava79} [Section 3.4,
equation (3.4.3)] provide the joint density of $H$ and $D=\operatorname
{diag}(d_i)$ in the standard Wishart case (which applies to $T$) as
%
\begin{eqnarray}\label{singdist1}
&&
f_{H,D}(H,D)
\nonumber\\[-8pt]\\[-8pt]
&&\qquad=C(p,n)|D|^{(n-p-1)/2} \biggl[\operatorname{etr} \biggl(-\frac12 D \biggr)
\biggr] \biggl[\prod_{i<j}(d_i-d_j)
\biggr]g_p(H)\nonumber
\end{eqnarray}
for constants $C(p,n)$ and functions $g_p$. Therefore, $H$ and $D$ are
independent. In the rank-deficient case ($p>n$), \citet{srivastava03}
(Section 3) provides an equivalent expression which, in
the singular Wishart case, gives
%
\begin{eqnarray}\label{singdist2}\qquad
&&
f_{H_1,D_1}(H_1,D_1)
\nonumber\\[-8pt]\\[-8pt]
&&\qquad=K(p,n)|D_1|^{(p-n-1)/2} \biggl[\operatorname{etr} \biggl(-\frac
12 D_1 \biggr) \biggr] \biggl[\prod_{i<j}(d_i-d_j)
\biggr]g_{n,p}(H_1)\nonumber
\end{eqnarray}
for constants $K(p,n)$ and functions $g_{n,p}$, so, again, we find
$H_1$ and $D_1$ independent by factorization. Now, $\lambda^+_{\min}$
is a function, in the full rank case (resp., rank-deficient case), of
only $D^{-1}$ (resp., $D_1^{-1}$), and we can write $T^+T=H'H$ (resp.,
$T^+T=H'_1H_1$), so $\lambda^+_{\min}$ and $T^+T$ are independent.
Being functions of $S$, they are also both independent of $X$. Now, the
nonzero eigenvalues of $T^+$ are the inverses of the nonzero
eigenvalues of $T$, a general fact about Moore--Penrose pseudo-inverses.
Therefore, denoting the largest eigenvalue of $T$ as
$\lambda_{\mathrm{max}}$, we can split up the expectations in
(\ref{mom4}) and get the bound
%
\begin{equation}\label{mom5}
\leq C^2_1\bigl\llvert 3+p-2 \operatorname{tr}
\bigl(SS^{+}\bigr)+n\bigr\rrvert E [ \lambda_{\mathrm{max}} ] E \biggl[
\frac{X'QX}{X'RX} \biggr] + 4C_1C_2.
\end{equation}

Now, it follows from positive semi-definiteness of $T$ that $
E [\lambda_{\mathrm{max}}  ]\leq
E [\operatorname{tr}(T) ]$. If $n\geq p$,
$\operatorname{tr}(T)\sim\chi^2_{pn}$ [cf. \citet{muirhead82},
Theorem 3.2.20] and so $E [\operatorname{tr}(T) ]=pn<\infty$.
If $p>n$, recall we can write $T=Z'Z$ for $Z\sim N_{n\times
p}(0,I_n\otimes I_p)$ by definition of the Wishart distribution; and
$ZZ'\sim\mathrm{Wishart}_n(p,I_n)$ so that
$\operatorname{tr}(T)=\operatorname{tr}(ZZ')\sim\chi^2_{pn}$; so,
again, $E [\operatorname{tr}(T) ]=pn<\infty$. Therefore, in
either case, $ E [\lambda_{\mathrm{max}}  ]\leq pn <
\infty$.

We still have to check that the expectation involving $X$, $Q$ and $R$
in (\ref{mom5}) is finite. Let $r=\operatorname{rk}(R)=\operatorname
{rk}(Q)=\operatorname{rk}(S)$ and write the spectral decomposition of
$(T^+TA)$ as $U\Lambda U'$, with $\Lambda=\operatorname
{diag}(L,0_{(p-r)})$ where $L$ is the vector of the $r$ nonzero
eigenvalues of $(T^+TA)$.
Then $R=(T^+TA)^+(T^+TA)=U\operatorname{diag}(I_r,0_{(p-r)}) U'$; let
us define the $p\times(p-r)$ matrix $E=U [0_{(p-r)\times r}
I_{(p-r)}]'$, that is,
so that $RE=0$ and $E$ has full column rank $p-r$. Notice that
$QE=AT^+TA U [0_{(p-r)\times r} I_{(p-r)}]'
=AU\Lambda U'U [0_{(p-r)\times r} I_{(p-r)}]'=0$. Since $Q$ and $R$
are symmetric positive semidefinite, we can use results in \citet
{magnus90} [Theorem 1(i) with
$A=Q$ and $B=R$] to conclude that
\[
E \biggl[ \frac{X'QX}{X'RX} \biggr] <\infty.
\]
This concludes the proof of the theorem.
\end{pf}

\section{Numerical study}\label{secnumstud}

This section provides some numerical results to showcase the
improvement in risk of the minimax estimator over the usual estimator.
More precisely, we compared the James--Stein estimator in (\ref{cor2})
given by
\[
\delta^{\mathrm{JS}}= \biggl(I-\frac{(n-2)SS^+}{(p-n+3)X'S^+X} \biggr)X
\]
and the usual estimator $\delta^0=X$ under invariant loss. (In
addition, we considered the positive James--Stein estimator to be
discussed in Section~\ref{seccomments}.)
The empirical approximations of the invariant risk of these estimators
were\vadjust{\goodbreak} plotted for $p=10, 20, 50$ and $n=\frac{p}{2}, p-1$. Three
covariance matrix structures were considered:\vspace*{9pt}

\textit{Spiked}: A diagonal matrix with the first $p/2$ diagonal elements
equal to~1, and the last $p/2$ equal to 10.

\textit{Autoregressive}: Autoregressive covariance matrices of the form
\[
\Sigma= %
\pmatrix{ 1 & \rho& \rho^2 &
\cr
\rho& 1 & \rho&
\cr
\rho^2 & \rho& 1 &
\cr
& & & \ddots }
\]
for $\rho= 0.5$.

\textit{Block diagonal}: Block diagonal matrices with $p/2$ blocks of the form
$ \bigl({1 \atop\rho}\enskip{\rho\atop1} \bigr)$
for $\rho=0.5$.\vspace*{9pt}

\noindent In all cases, the true mean was chosen as $\theta\propto(1,\ldots,1)$.

\begin{figure}

\includegraphics{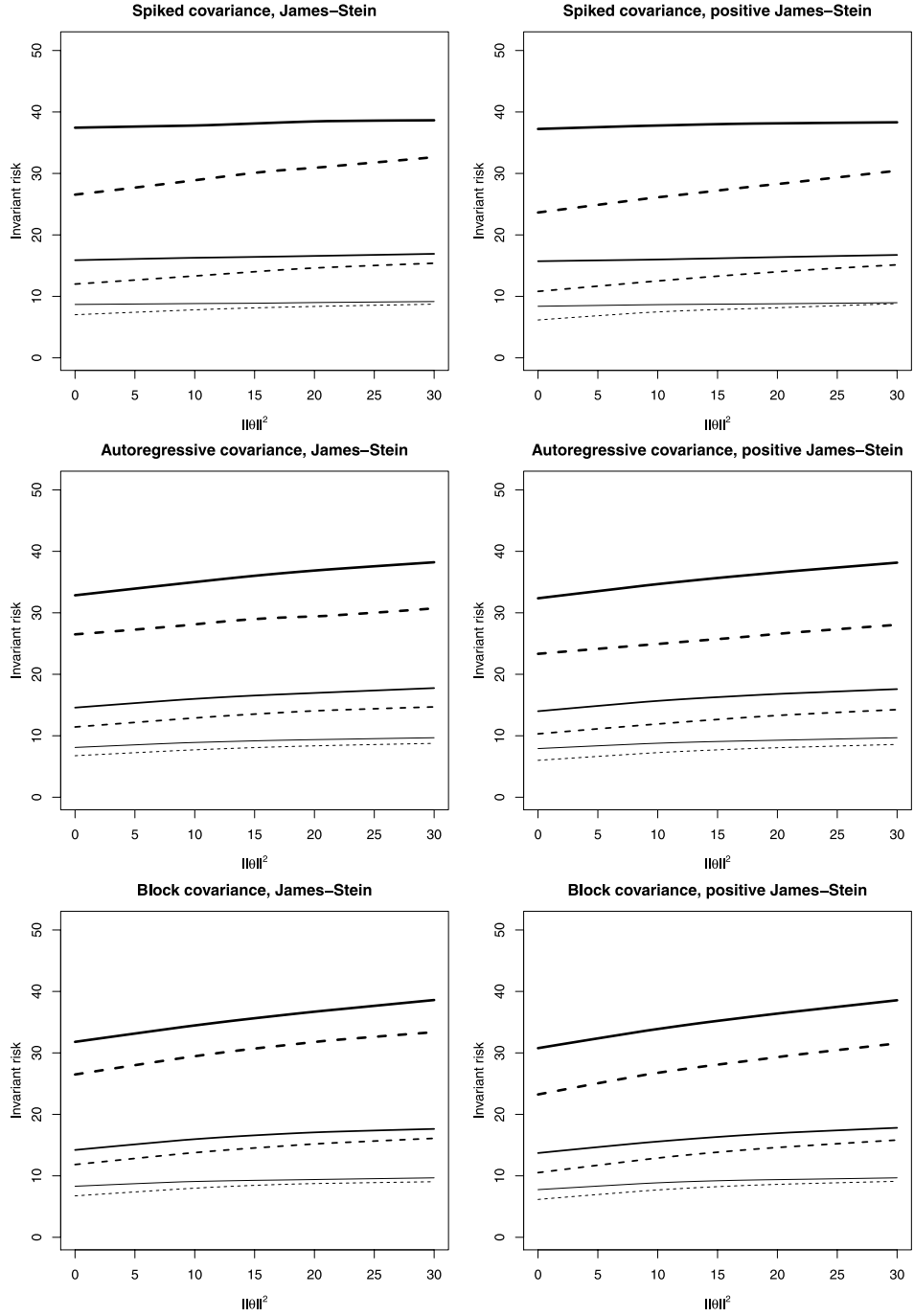}

\caption{The risk function plots of $\delta_a^{\mathrm{JS}}$ and
$\delta_a^{\mathrm{JS}+}$ for $a=(n-2)/(p-n+3)$ are in the left and
right columns, respectively. The lines, from thinnest to thickest,
are for $p=10, 20$ and $50$. The solid and dashed lines are,
respectively, for $n=p/2$ and $n=p-1$.}
\label{simplots}
\end{figure}

We remind the reader that the risk of the trivial estimator is always
$p$, regardless of $\theta$ or $\Sigma$. With this in mind, we see
from Figure~\ref{simplots} that in all six scenarios the pattern of
%
domination of the new estimator is similar to one of the usual
James--Stein estimators. Also note that, as predicted by the theoretical
results, the domination decreases as the smaller $n$ tends to $p$.

\section{Comments}\label{seccomments}

An interesting property of the Moore--Penrose inverse is that for any
$A$, $AA^+$ is the matrix that projects onto the subspace spanned by
$A$ (its column space). It follows that the proposed generalized
Baranchik estimator can be expressed as
%
\begin{eqnarray}\label{projform}
\delta_r(X,S)&=&\bigl(I-SS^+\bigr)X+ \biggl(1-\frac{r(X'S^+X)}{X'S^+X}
\biggr)SS^+X
\nonumber\\[-8pt]\\[-8pt]
&=&P_{S^\perp}X+ \biggl(1-\frac{r(X'S^+X)}{X'S^+X} \biggr)P_S
X,\nonumber
\end{eqnarray}
where $P_S=SS^+$ and $P_{S^\perp}=I-SS^+$ are the projection matrices
onto the column space of $S$ and its orthogonal complement,
respectively. In terms of the kernel and image of the symmetric matrix
$S$, $\operatorname{Ker}(P_{S^\perp})= \operatorname{Im}(S)$ and
$\operatorname{Im}(P_{S^\perp})=\operatorname{Ker}(S^+)$. When
$p>n$, this means we can
interpret our estimator as applying shrinkage only on the component of
$X$ in the subspace spanned by our covariance matrix estimator $S$. In
particular, note that the estimator
$P_S\delta_r(X,S)= (1-\frac{r(X'S^+X)}{X'S^+X} )P_S X$
dominates $P_S X$ under invariant loss function (\ref{eq2}), since
$R(P_S\delta_r, \theta)- R(P_S X, \theta) = R(\delta_r, \theta)- R(X,
\theta)\ge0$ if $r$ satisfies the conditions of
Theorem~\ref{domination}. This suggests there might be an easier, more
abstract proof of Theorem~\ref{domination}, one not relying on brute
computations but on the already known full rank $S$ case, although we
have not been able to obtain such a result.\vadjust{\goodbreak}

A natural extension of the James--Stein estimator, $\delta^{\mathrm{JS}}_a$ in
(\ref{cor2}), is a positive-part-type James--Stein estimator. The form
of the estimator in (\ref{projform}) suggests
%
\begin{equation}
\label{ppjs} \delta_a^{\mathrm{JS}+}= \bigl(I-SS^+ \bigr)X+
\biggl(1-\frac
{a}{X'S^+X} \biggr)_{+}SS^+X,
\end{equation}
where $b_+ = \max(b,0)$. Simulation evidence from Figure \ref
{simplots} suggests that for $a=(n-2)/(p-n+3)$, $\delta_a^{\mathrm
{JS}+}$ dominates $\delta^{\mathrm{JS}}_a$ under invariant loss.

One of the interesting differences between the $n>p$ and $p>n$ cases is
the reversal of the roles of $p$ and $n$. This is essentially due to
the distribution of the singular values of $S$. Recall that for
$S=ATA$, $T\sim W_p(n,I_n)$. We can write the singular value
decomposition of $T$ as $T=H'DH$, but we can also write it as $T=H'_1D_1H_1$,
where $H_1$ is semi-orthogonal ($H_1H'_1=I)$ and $D_1$ is the matrix of
the positive eigenvalues of $T$. If $T$ has full rank (i.e., $n\geq
p$), this coincides with the singular value decomposition of $T$. In
the full rank case the joint density of $H$ and $D$ is given in (\ref
{singdist1}), whereas in the rank-deficient case ($p>n$) joint density
is given by (\ref{singdist2}), from which stems the reversal of the
roles of $p$ and $n$.

In the heteroscedastic normal mean estimation problem, \citet
{JamesStein61proBer} used the loss function that was weighted by the
inverse of the variances and, consequently, the problem is essentially
transformed to the homoscedastic case under ordinary squared error
loss. Similarly, in this article, we used the invariant loss function
in (\ref{eq2}), therefore skirting a somewhat subtle issue. In the
heteroscedastic setting where there are differing coordinate variances,
minimax estimation and Bayes (or empirical Bayes) estimates can be
qualitatively different. It turns out that minimax estimators in
general shrink most on the coordinates with smaller variances, while
Bayes estimators shrink most on large variance coordinates. \citet
{Brown75} shows that the James--Stein shrinkage estimator does not
dominate the $X$ when the largest variance is larger than the sum of
the rest. Moreover, \citet{Casella80} points out that the
James--Stein shrinkage estimator may not be a desirable shrinkage
estimator under heteroscedasticity even when it is minimax. \citet
{Morris09} and \citet{BNX12} give an excellent perspective on
minimaxity of the shrinkage estimator from Bayes and empirical Bayes
points of view. Consequently, it would be of interest to examine the
shrinkage patterns of the proposed estimates in the case of a
noninvariant loss function and assess how well the invariant loss works
for $p>n$ applications.

One can imagine an extension of the results of this article beyond the
normal distribution setting. Consider a model with the joint
density for $(X, S)$ the form
%
\begin{equation}
\label{ess} f \bigl(\operatorname{tr} \Sigma^{-1}\bigl[(X - \theta) (X -
\theta)^{\prime} + S\bigr] \bigr),
\end{equation}
where the $p \times1$ location vector $\theta$ and the $p \times p$
scale matrix
$\Sigma$ are unknown. In the setting of $p\leq n$, \citet
{fourdrinier03} and
\citet{KubokawaSrivastava01} give some results on improved
location estimation for elliptically symmetric distributions. For more
on elliptical symmetry and the various
choices of $f(\cdot)$ in (\ref{ess}), see \citet
{FangKotzNg90book}; the class in
(\ref{ess}) contains models such as the multivariate normal, $t$- and
Kotz-type distributions.

Finally, simulation study reveals that, when $p$ is much larger than
$n$, the estimate of $\Sigma$ and $\Sigma^{-1}$ are quite poor. This
observation agrees with \citet{KubSriv08}, where \citet
{Haff79JMA}-type improved estimates of $\Sigma$ are proposed. It would
be of interest to use an improved estimator of $\Sigma$ in
$\delta_r(X,S)$ in (\ref{deltar}). As pointed out in the testing
context by \citet{SrivFuj06} and \citet{Sriv07}, a
shortcoming of $S^+$ is that the associated estimator is only
orthogonally invariant, while the sample mean vector is invariant.

\section*{Acknowledgments}

The authors are grateful to the Associate Editor and referees for
helpful comments that strengthened the exposition and scope of this
paper.



\printaddresses

\end{document}